\documentclass[11pt,a4paper]{article}
\usepackage{fullpage}
\usepackage{xspace,xcolor}
\usepackage{setspace}
\usepackage{afterpage}
\usepackage{graphicx,epsfig,epstopdf}
\usepackage{amscd}
\usepackage{ifthen}
\usepackage{gensymb}
\usepackage{amsmath}
\usepackage{amsfonts}
\usepackage{fancyhdr}
\usepackage{amssymb}
\usepackage{dcolumn}
\usepackage{mathtools}
\usepackage{bm}
\usepackage{array}
\usepackage{cite}
\usepackage{caption}
\usepackage{subcaption}
\usepackage{multirow}
\usepackage{nicematrix}
\usepackage{blindtext}
\usepackage{hyperref}
\usepackage{nameref}
\usepackage{csquotes}
\usepackage{longtable}

\newcommand\blfootnote[1]{%
	\begingroup
	\renewcommand\thefootnote{}\footnote{#1}%
	\addtocounter{footnote}{-1}%
	\endgroup
}

\begin{document}
	\allowdisplaybreaks
	
	\thispagestyle{empty}
	
\begin{center}
\doublespacing
\medskip
\medskip
{\LARGE \bfseries Hybrid GFD-RBF Method for Convection-Diffusion Problems}\\
\singlespacing
\vspace*{0.25cm}
{\large  Priyal Garg$^\dagger$ and T.V.S. Sekhar*\blfootnote{*Corresponding Author}}\\
\vspace*{0.25cm}
{School of Basic Sciences, Indian Institute of Technology Bhubaneswar, Odisha-752050, India}\\
$^\dagger$\it{a21ma09003@iitbbs.ac.in, *sekhartvs@iitbbs.ac.in}

\end{center}

\section*{\Large{Abstract}}
\begin{singlespace}
In this paper, we present a meshless hybrid method combining the Generalized Finite Difference (GFD) and Finite Difference based Radial Basis Function (RBF-FD) approaches to solve non-homogeneous partial differential equations (PDEs) involving both lower and higher order derivatives. The proposed method eliminates the need for mesh generation by leveraging the strengths of both GFD and RBF-FD techniques. The GFD method is robust and stable, effectively handling ill-conditioned systems, while the RBF-FD method excels in extending to higher-order derivatives and higher-dimensional problems. Despite their individual advantages, each method has its limitations. To address these, we developed a hybrid GFD-RBF approach that combines their strengths. Specifically, the GFD method is employed to approximate lower order terms (convective terms), and the RBF method is used for higher order terms (diffusive terms). The performance of the proposed hybrid method is tested on both linear and nonlinear PDEs, considering uniform and non-uniform distributions of nodes within the domain. This approach demonstrates the versatility and effectiveness of the hybrid GFD-RBF method in solving second and higher order convection-diffusion problems.\\

Keywords : {\it{Hybrid Meshless Method, Generalized Finite Difference Method, Radial Basis Function, GFD-RBF Method}}.
\end{singlespace}

\section{Introduction}
Non-homogeneous PDEs involving lower and higher order derivatives on complex or moving boundaries are frequently used to model various physical phenomena. Addressing the challenges posed by complex boundaries, meshless methods offer greater flexibility compared to mesh-based methods. By discretizing derivatives at scattered points without connecting them with edges, meshless methods are well-suited for problems involving moving boundaries, such as fluid-structure interaction. For boundary layer flows, adaptive local refinement is also straightforward, as there is no need to account for grid properties like control volume aspect ratio or element skewness. Recently, meshless numerical manifold method is also used for various physical problems\cite{1,2}. Localized meshless methods, such as GFD and RBF-FD have proven effective for these tasks.\\

GFD method is based on the Taylor series expansion and the moving least squares (MLS) technique for approximating the derivatives of an unknown function. Theoretical advancements in the GFD method were made by Jensen \cite{GFD theory 1}, Perrone and Kao \cite{GFD theory 2}, and Benito et al. \cite{GFD theory 3}, with subsequent analyses by Benito et al. \cite{GFD theory 4} and Zheng and Li \cite{GFD theory 6}. Applications of the GFD method have been demonstrated by Mochnacki and Majchrzak \cite{GFD app 2} and Urena et al. \cite{GFD app 3} across various fields. Gavete et al. \cite{GFD theory 5} applied the method to solve nonlinear elliptic PDEs, while Urena et al. \cite{para, hyper} tackled nonlinear parabolic and hyperbolic PDEs. GFD method has notable strengths, including its robustness and stability in handling ill-conditioned systems \cite{GFD app 4,GFD app 5}. However, it requires the computation of all derivatives of a particular order, even when only a few is needed for a specific PDE. Additionally, improving accuracy by including higher order terms in the Taylor series becomes computationally expensive due to the evaluation of numerous higher order derivatives. Consequently, GFD method is less efficient for solving higher order PDEs or achieving higher order accuracy. RBF-FD method \cite{RBF theory 1, RBF theory 2, RBF theory 3, RBF theory 4}, on the other hand, uses radial basis function interpolation to approximate derivatives at given points. It is well-suited for higher-order PDEs and high-dimensional problems, as demonstrated in various applications \cite{RBF app 1, RBF app 2} and modifications to the classical method \cite{RBF app 3, RBF app 4, RBF app 5}. The key advantage of the RBF-FD method lies in its extensibility and accuracy. However, its performance depends on selecting an optimal shape parameter. Further, if the number of nodes is increased beyond a certain limit, it can lead to discretization errors reaching a constant value, which is known as the stagnation or saturation phenomenon \cite{RBF app 4}.\\

 In general, numerical solutions for higher order PDEs are obtained by rewriting them as systems of lower order PDEs. However, this approach increases the number of unknowns and poses challenges in handling boundary conditions, particularly for complex boundaries. To address these issues, we have developed a hybrid GFD-RBF method that directly solves higher order PDEs by retaining the strengths of both GFD and RBF methods. In this hybrid approach, lower order derivatives (e.g., convection terms) are approximated using the GFD method, while higher order derivatives (e.g., diffusion terms) are computed using the RBF-FD method. The method has been tested on both linear and nonlinear PDEs, with a comparative analysis conducted with the standard second-order central finite difference (CD2) method for uniformly distributed nodes and GFD method for both uniform and non-uniform distributions.

\section{Numerical Approximation}\label{Numerical Approximation}
Consider a second-order non-homogenous partial differential equation of Dirichlet type:
\begin{equation}\label{PDE}
	\mathcal{L}(U(\textbf{x}))=f(\textbf x)~~~\rm in~~ \Omega
\end{equation}
with boundary conditions
\begin{equation}
	U(\textbf{x})=g(\textbf x)~~\rm on~~\Gamma
\end{equation}
where, $\textbf x=(x,y)$.\\

\begin{figure}[h!]
	\centering
	\includegraphics[width=8cm,height=8cm,angle=0]{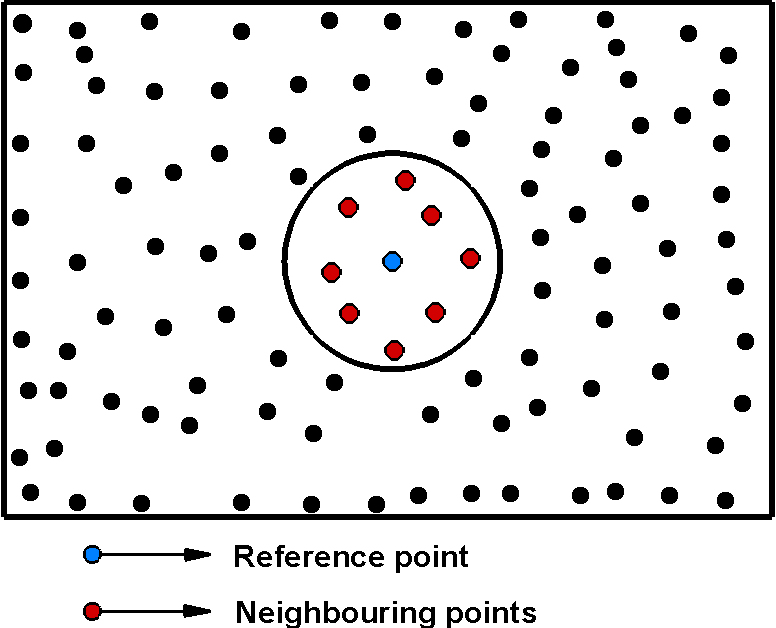}
	\caption{Domain showing star with central node and neighbouring nodes}
	\label{fig:star}
\end{figure}

Let $n$ be the total number of nodes in the domain. The motive is to find the value of $U$ at the internal nodes. For this, at each internal node, say $\textbf x_0$, define a composition with $\textbf x_0$ as the central node and $\textbf x_1,\textbf x_2, \dots,\textbf x_{N-1}$ as its neighbouring nodes. This composition is termed as star which is shown in Fig. \ref{fig:star}. The star has been selected as per the quadrant rule \cite{selection of points}.\\

In order to find the approximate value of $U(\textbf x)$, say $u(\textbf x)$, given PDE has to be discretized to get the corresponding difference equation. In hybrid GFD-RBF method, we propose that the lower order derivatives are approximated by GFD and the higher order derivatives by RBF-FD. For example, if a second order PDE is considered, the first order derivatives are approximated by GFD while the second order derivatives by RBF. Similarly, if a third order PDE is considered, first order derivatives are approximated by GFD and higher ones by RBF or derivatives up to second order are approximated by GFD and higher ones by RBF.

\subsection{Approximation of the First Order Derivatives by GFD}
For a particular star, let $U(\textbf x_0)$ is the value of the function at the central node and  $U(\textbf x_i),~i=1,\dots,N-1$ are the function values at the neighbouring nodes, then according to the Taylor series expansion, we have:
\begin{equation}
U_i=U_0+h_i\frac{\partial U_0}{\partial x}+k_i\frac{\partial U_0}{\partial y}+\dots,~~~~i=1,\dots,N-1,
\end{equation}
$h_i=x_i-x_0,~k_i=y_i-y_0$. If $u(\textbf x)$ be the approximate value of $U(\textbf x)$, then to get the approximations of first order, we truncate the series after first order terms. It can be represented in matrix form as:

\begin{equation}\nonumber 
	\begin{bmatrix}
		h_1 & k_1 \\
		h_2 & k_2 \\
		\vdots & \vdots\\
		h_{N-1} & k_{N-1} 
	\end{bmatrix}
	\begin{bmatrix}
		\frac{\partial u_0}{\partial x}\\
		\frac{\partial u_0}{\partial y}\\
	
	\end{bmatrix}=
	\begin{bmatrix}
		u_1-u_0\\
		u_2-u_0\\
		\vdots\\
		u_{N-1}-u_0
	\end{bmatrix},
\end{equation}
or 
\begin{equation}\label{gfd_matrix}
	AD_{\textbf u_0}=\Delta \textbf u.
\end{equation}
Let
\begin{equation}
	\textbf R=\Delta \textbf u-AD_{\textbf u_0}
\end{equation}
denotes the residue of Eq. (\ref{gfd_matrix}). The residuals can be minimized by applying the MLS method. If $W=diag(w_1,w_2, \dots ,w_{N-1})$ be the matrix of the weights, then according to the method we have to solve the following system of equations:
\begin{equation}\label{Derivative of residue}
	\frac{\partial \textbf R^{'} W^{2}\textbf R}{\partial D_{\textbf u_0}}=0.
\end{equation}
The approximate values of the derivative can be obtained by solving Eq. (\ref{Derivative of residue}) as:
\begin{equation}
	D_{\textbf u_0}=B\Delta \textbf u
\end{equation}
where,\\
$B=S^{-1}A^{'}W^{2}$, $S=A^{'}W^{2}A$. Using the quadrant criteria, the neighboring nodes around the central node are chosen in such a way that the matrix $A$ contains a base of $\mathbb{R}^{2}$, and hence $S$ is invertible.\\

Let $\textbf u=(u_0, u_1,\dots, u_{N-1})$ and
$\bar{B}=\begin{bmatrix}
	-\sum_{i=1}^{N-1}b_{1,i} & b_{1,1} & b_{1,2} & \cdots & b_{1,N-1} \\
	-\sum_{i=1}^{N-1}b_{2,i} & b_{2,1} & b_{2,2} & \cdots & b_{2,N-1} 
\end{bmatrix}$\\
where, $b_{i,j}$ are the elements of $B$. Now the expression for the derivatives will be given by
\begin{equation}\label{convection}
	D_{\textbf u_0}=B\Delta \textbf u=\bar{B}\textbf u=\begin{bmatrix}
		-u_0\sum_{i=1}^{N-1}b_{1,i}+\sum_{i=1}^{N-1}b_{1,i}u_{i}\\
		-u_0\sum_{i=1}^{N-1}b_{2,i}+\sum_{i=1}^{N-1}b_{2,i}u_{i}
	\end{bmatrix}.
\end{equation}

The above equation will give the discretized form of the first order derivatives at an internal node $(\textbf x_0)$ of a star of the domain. The same process is repeated by considering each internal node as the central node and assigning a star to that node.
\subsection{Approximation of the Second Order Derivatives by RBF-FD}
Let's consider the same star as for GFD with $N$ distinct nodes. If $u(\textbf x)$ be the approximate value of $U(\textbf x)$, then the RBF interpolation is given by
\begin{equation}\label{rbf_equation}  
u(\textbf x) \approx s(\textbf x)=\sum_{i=0}^{N-1}\lambda_i\phi(\parallel\textbf x-\textbf x_i\parallel,\epsilon)+\sum_{j=0}^{M-1}\mu_jp_j(\textbf x)
\end{equation}
where, $\phi$ is the RBF function with $\epsilon$ as the shape parameter and $\{p_j(\textbf x)\}^{M-1}_{j=0}$ is a basis for $\Pi_M(\mathbb{R}^2)$ (space of all 2-variate polynomials of degree less than or equals to $M$) which is required for the interpolation problem to be well-posed\cite{Rbf extra}. Putting $u(\textbf x_i)=s(\textbf x_i),~i=0,\dots,~N-1$ in Eq. (\ref{rbf_equation}), we get $N$ linear equations and therefore $M$ extra conditions are required to solve the equation. These extra conditions are chosen by taking the expansion coefficient vector $\lambda\in\mathbb{R}^N$ orthogonal to $\Pi_M(\mathbb{R}^2)$, i.e.,
$$\sum_{i=0}^{N-1}\lambda_ip_j(\textbf x_i)=0,~~~~~j=0, \dots, M-1.$$
This will lead to a system of linear equations which can be represented as:
\begin{equation}\label{rbf_matrix}
\begin{bmatrix}
\bm{\phi} & \textbf{p}\\
\textbf{p} & \textbf{0}
\end{bmatrix}
\begin{bmatrix}
\bm{\lambda}\\
\bm{\mu}
\end{bmatrix}=
\begin{bmatrix}
	\textbf{u}\\
	\textbf{0}
\end{bmatrix}
\end{equation}
where, $\bm{\phi}\coloneq\phi(\parallel\textbf x_j-\textbf x_i\parallel,\epsilon),~i,j=0,\dots,N-1$ and $\textbf{p}\coloneq p_{j}(\textbf x_{i}),~j=0,\dots,M-1,~i=0,\dots,N-1$.\\

In order to derive the derivatives of the function $u(\textbf x)$ at a given point, Lagrange form of RBF interpolant is used which is given by
\begin{equation}\label{lagrange_form}
u(\textbf x) \approx s(\textbf x)=\sum_{i=0}^{N-1}\psi_{i}(\textbf x)u(\textbf x_i)
\end{equation}
where, $\psi_{i}(\textbf x)$ satisfies the condition of Kronecker delta, i.e.,
\begin{equation}
\psi_{i}(\textbf x_j)=\delta_{ij},~~~j=0,\dots,N-1.
\end{equation}
The closed form representation for $\psi_i(\textbf x)$ can be obtained by considering that the
right-hand side vector of Eq. (\ref{rbf_matrix}) stems from each $\psi_i$'s. Hence,
\begin{equation}
\psi_i(\textbf x)=\frac{det(Q_i(\textbf x))}{det(Q)}
\end{equation}
where,
\begin{equation}\label{Qmatrix}
	Q=\begin{bmatrix}
	\bm{\phi} & \textbf{p}\\
	\textbf{p} & \textbf{0}
\end{bmatrix}
\end{equation}
and $Q_i(\textbf x)$ is same as the matrix $Q$, except that the $i$th row is replaced by the vector 
\begin{equation}\label{B_vector}
B(\textbf x)=[\phi(\parallel\textbf x-\textbf x_0\parallel,\epsilon),\dots,\phi(\parallel\textbf x-\textbf x_{N-1}\parallel,\epsilon)~|~p_0(\textbf x),\dots,p_{M-1}(\textbf x)].
\end{equation}

The representations through Eqs. (\ref{lagrange_form}-\ref{Qmatrix}) can be used to approximate the derivative of a function $\textbf Lu$ at a given point $\textbf x_j$, where $\textbf L$ is a second order linear differential operator. To derive the formulae for $\textbf L(u(\textbf x_0))$, a set of neighboring node of $\textbf x_0$ (say $\{\textbf x_0,\textbf x_1,\dots,\textbf x_{N-1}\}$) is considered. Let $\textbf L(u(\textbf x_0))$ be represented as a linear combination of the value of $u$ at $\textbf x_0$ and its neighbouring nodes, i.e.,
\begin{equation}\label{linear_form}
\textbf L(u(\textbf x_0))=\sum_{i=0}^{N-1}c_iu(\textbf x_i),~\text{for each}~\textbf x_0\in\Omega
\end{equation}
where, $c_i$'s are the weights to be calculated. So, after applying the operator $\textbf L$ to the Lagrange form of RBF interpolant in Eq. (\ref{lagrange_form}) for $\textbf x_0$, we get 
\begin{equation}\label{final_form}
\textbf L(u(\textbf x_0)) \approx \textbf L(s(\textbf x_0))=\sum_{i=0}^{N-1}\textbf L(\psi_{i}(\textbf x_0))u(\textbf x_i).
\end{equation}
On comparing Eqs. (\ref{linear_form}-\ref {final_form}), we get $c_i$'s as:
\begin{equation}
c_i=\textbf L(\psi_{i}(\textbf x_0)),~i=0,\dots,N-1.
\end{equation}
The weights can be computed by solving the following linear system:
\begin{equation}
Q[\textbf c/\bm{\mu}]=(\textbf L(B(\textbf x_0)))^T,~\text{for each}~\textbf x_0\in\Omega
\end{equation}
where, $Q$ and $B$ are given by Eq. (\ref{Qmatrix}) and Eq. (\ref{B_vector}) respectively and $c=[c_0,c_2,\dots,c_{N-1}]^T$ denotes the vector of weights.\\

Putting the values of weights in Eq. (\ref{linear_form}) will give the discretized form of the second order derivatives at an internal node of the domain. The same process is repeated for each internal node of the domain.\\

Substituting the value of first order derivatives obtained from GFD and the second order derivatives obtained from RBF-FD for a given node in Eq. (\ref{PDE}), we get an algebraic equation corresponding to that node. The same process is repeated for each internal node of the domain to get a system of equations in which the number of equations and unknowns will be equal to the number of internal nodes. This algebraic system is solved using the stabilized bi-conjugate gradient (Bi-CGSTAB) method to obtain the approximate value of the function $U$ at each node of the domain.

\section{Results and Discussion}	 
\begin{figure}[h]
	\begin{subfigure}{.5\textwidth}
		\centering
		\includegraphics[width=.95\linewidth]{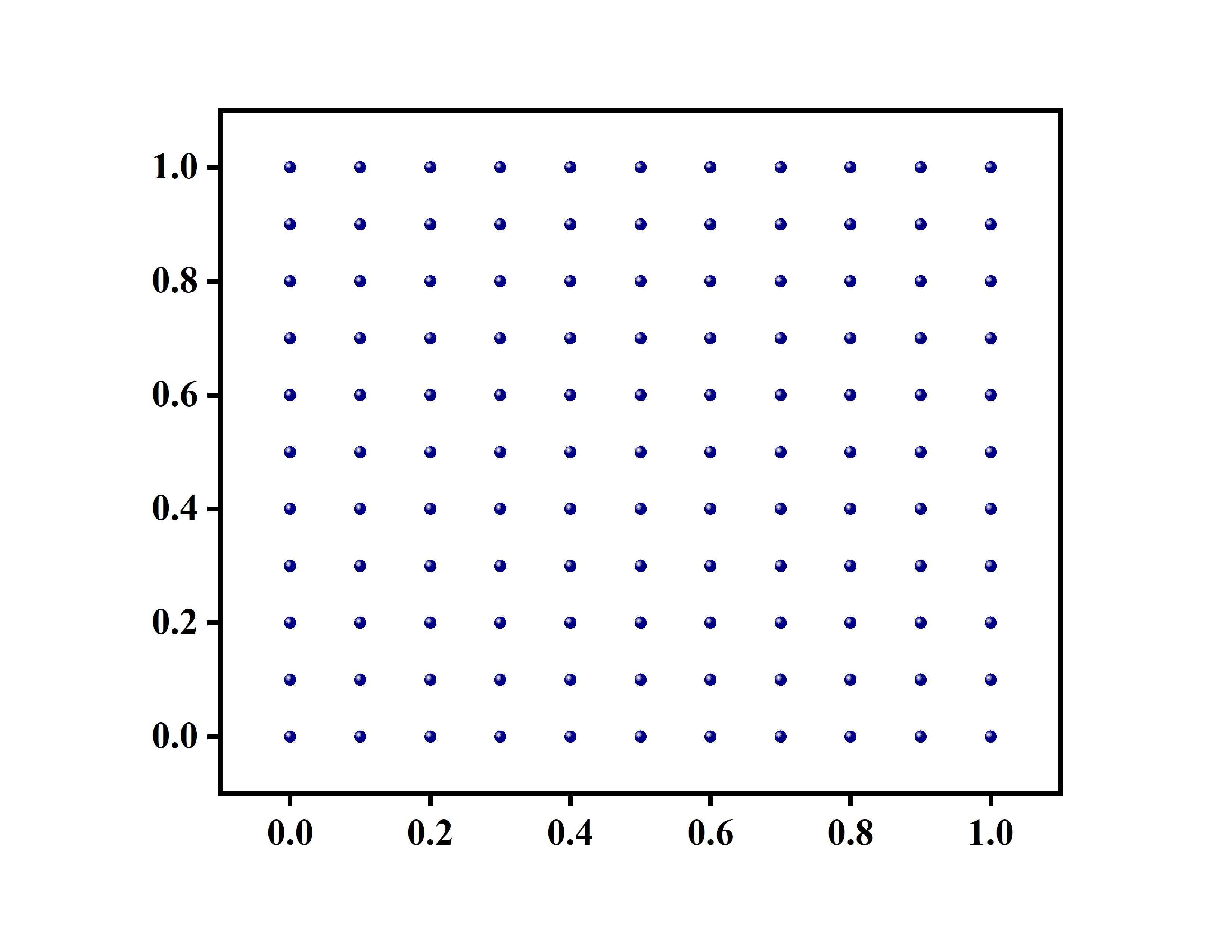}
		\caption{Uniform Nodes}
	\end{subfigure}%
	\begin{subfigure}{.5\textwidth}
		\centering
		\includegraphics[width=.95\linewidth]{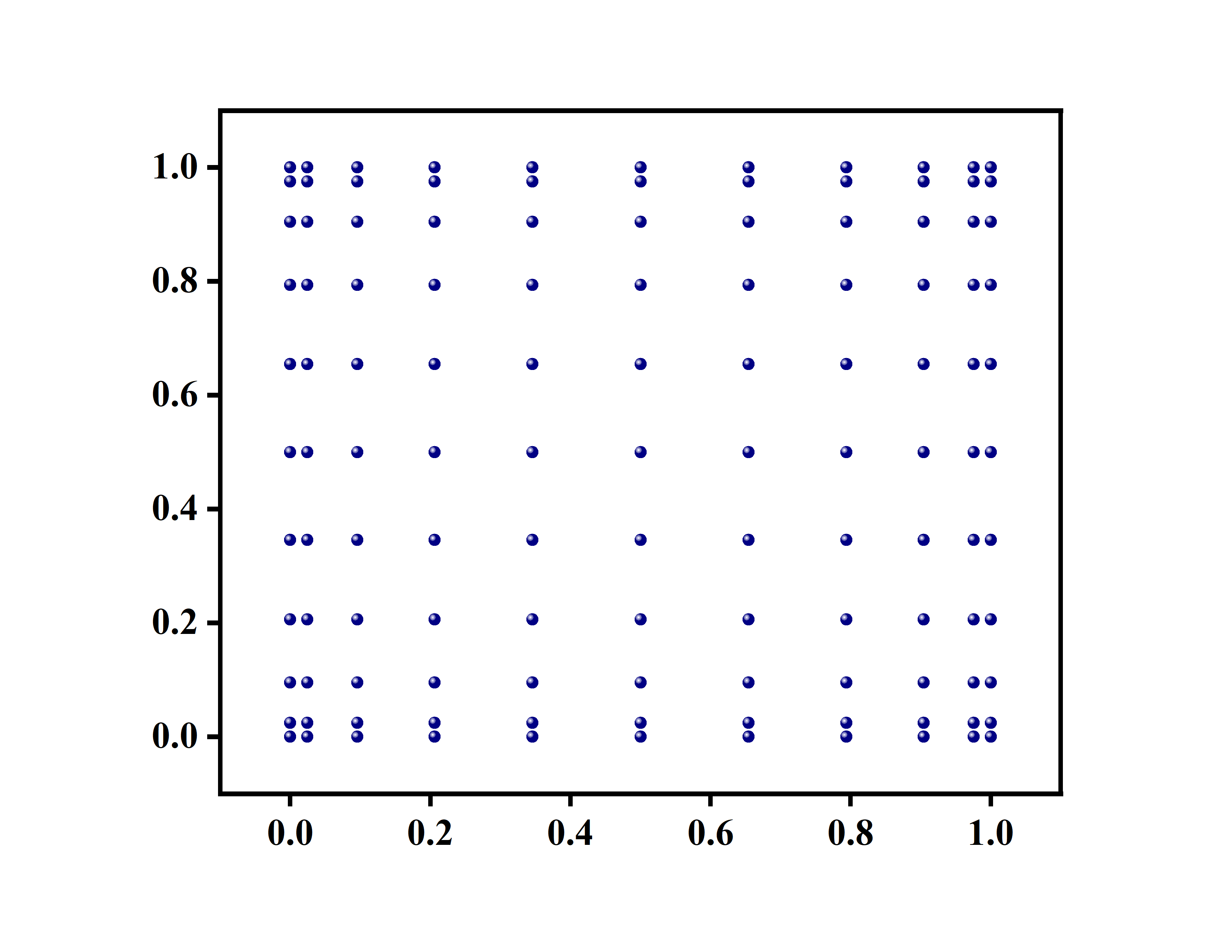}
		\caption{Chebyshev Nodes}
	\end{subfigure}
	\caption{Discretization of the Domain with 121 Nodes}
	\label{fig:Domain}
\end{figure}

The hybrid method has been applied to both linear and nonlinear PDEs with known analytical solutions. Boundary conditions and the source terms can be calculated from the analytical solution wherever required. Domains are discretized by considering uniform and Chebyshev distribution of nodes. The sample of 2D domains with chosen points is shown in Fig. \ref{fig:Domain}. Chebyshev nodes are generated by using the formula: 
\begin{equation}\nonumber
	x_i=\frac{1}{2}\left[1-\cos\left(\frac{i-1}{n_x-1}.\pi\right)\right],~~i=1,\dots,n_x,
\end{equation}
\begin{equation}\nonumber
	y_j=\frac{1}{2}\left[1-\cos\left(\frac{j-1}{n_y-1}.\pi\right)\right],~~j=1,\dots,n_y.
\end{equation}
The weight function for GFD is taken as $1/(distance)^2$ and multiquadric function ($\sqrt{1+\epsilon^2r^2}$)\cite{Multiquadric} is used for RBF-FD. 

\subsection{Illustrative Examples}
\textbf{Example 1: Linear Convection-Diffusion Equation}
$$-2\frac{\partial U}{\partial x}+2\frac{\partial U}{\partial y}-\frac{\partial^2 U}{\partial x^2} -\frac{\partial^2 U}{\partial y^2}=f(x,y),~~~~~~~~0\leq x,y\leq 1.$$
Exact solution:
$$U(x,y)=\frac{e^{2(1-x)}+e^{2y}-2}{e-1}.$$

\begin{table}[h!]
\caption{Comparison of Error, Order of Convergence and Condition Number of the Stiffness Matrix for Example 1 with $\epsilon=0.5$}
\vspace{-0.3 cm}
\begin{center}
		\vspace{-0.4 cm}
		\scriptsize
		\begin{tabular}{c c c c c c c c c}
			\hline
			\rule{0pt}{12pt} \multirow{2}{*}{Nodes} & \multicolumn{3}{c}{Hybrid GFD-RBF}  &\multicolumn{3}{c}{CD2} &\multicolumn{2}{c}{Ref.\cite{Comparision}*(RBf-FD)}\\
			\cline{2-9}
			\rule{0pt}{12pt} & Error & Order & C.N. & Error & Order & C.N. & Error & Order \\
			\hline
			\rule{0pt}{12pt} $121$ & $3.4543\times 10^{-5}$ & - & $38.3205$ & $2.1827\times 10^{-4}$ & - & $38.3125$ & $1.3257\times 10^{-3}$ & - \\
		\rule{0pt}{12pt}$ 441$ & $8.9838\times 10^{-6}$ & $1.9429$ & $154.8495$ & $5.8404\times 10^{-5}$ & $1.9019$ & $154.8409$ & $3.7435\times 10^{-4}$ & $1.82$\\
		\rule{0pt}{12pt}$ 1681$ & $2.2998\times 10^{-6}$ & $1.9658$ & $620.9606$ & $1.9312\times10^{-5}$ & $1.5965$ & $620.9518$ & $9.1949\times 10^{-5}$ & $2.02$\\
		\hline
	\end{tabular}
		\label{P1 CD2}
	\end{center}
\vspace{-0.3cm}
\scriptsize{*Number of nodes taken by the authors for RBF-FD are $81$, $289$ and $1089$.}
\end{table}

In Example 1, the proposed hybrid method has been applied to linear convection-diffusion equation by discretizing the convective terms using GFD and the diffusive terms by RBF with five nodes per star (NG = 5, NR = 5) as described in Section \ref{Numerical Approximation}. To analyse the efficiency of the hybrid method, the root mean square (RMS) error (using analytical solution), order of convergence of the method and the condition number (C.N.) of the stiffness matrix with uniform distribution of nodes are calculated  and presented in Table \ref{P1 CD2}. The order of convergence of the method is given by
\begin{equation}\nonumber
	q\approx\frac{\log(e_{\text{old}}-e_{\text{new}})}{\log(h_{\text{old}}-h_{\text{new}})}
\end{equation}
where, $e_{\text{old}}$ and $e_{\text{new}}$ are the RMS errors with nodal spacing $h_{\text{old}}$ and $h_{\text{new}}$ respectively. The condition number of stiffness matrix is calculated since the stability of a system of algebraic equations is closely related to the condition number of its coefficient matrix. The proposed hybrid method is compared with CD2 and RBF-FD\cite{Comparision} methods and the results are given in Table \ref{P1 CD2}. It is found that the hybrid method gives more accurate results and it is second order accurate. In fact, the developed hybrid method is giving more satisfactory results with less number of nodes (e.g. nodes-$121$, error- $3.4543\times 10^{-5}$) in the domain when compared with RBF-FD (e.g. nodes-$289$, error- $3.7435\times 10^{-4}$). The condition number of the stiffness matrix of the hybrid method and CD2 is also comparable with each other. \\

\begin{table}[h!]
	\caption{Comparison of Error and Order of Convergence  for Different Nodes per Star for Uniform Distribution of Nodes for Example 1 with $\epsilon=0.5$}
	\vspace{-0.3 cm}
	\begin{center}
		\vspace{-0.4 cm}
		\begin{tabular}{c c c c c c c}
			\hline
			\rule{0pt}{12pt}\multirow{2}{*}{Nodes} &\multicolumn{2}{c}{$N_G=5,~N_R=5$} &\multicolumn{2}{c}{$N_G=9,~N_R=5$}& \multicolumn{2}{c}{$N_G=9,~N_R=9$}\\
			\cline{2-7}
			\rule{0pt}{12pt} & Error & Order & Error & Order & Error & Order\\
			\hline
			\rule{0pt}{12pt}$ 121$ & $3.4543\times 10^{-5}$ & - & $3.4601\times 10^{-5}$ & - &  $1.4893\times 10^{-4}$ & - \\
			\rule{0pt}{12pt}$ 441$ & $8.9838\times 10^{-6}$ & $1.9429$ & $8.9878\times 10^{-6}$ & $1.9447$ & $4.1394\times 10^{-5}$ & $1.8471$ \\
			\rule{0pt}{12pt}$ 1681$ & $2.2998\times 10^{-6}$ & $1.9658$ & $1.2300\times 10^{-6}$ & $1.9662$ & $2.4561\times 10^{-5}$ & $0.7530$ \\
			\hline
		\end{tabular}
		\label{P1 Error Uniform}
	\end{center}
\end{table}

\begin{table}[h!]
	\caption{Comparison of Error and Order of Convergence for Different Nodes per Star for Chebyshev Distribution of Nodes for Example 1 with $\epsilon=0.9$}
	\vspace{-0.3 cm}
	\begin{center}
		\vspace{-0.4 cm}
		\begin{tabular}{c c c c c c c}
			\hline
			\rule{0pt}{12pt}\multirow{2}{*}{Nodes} &\multicolumn{2}{c}{$N_G=5,~N_R=5$} &\multicolumn{2}{c}{$N_G=9,~N_R=5$}& \multicolumn{2}{c}{$N_G=9,~N_R=9$}\\
			\cline{2-7}
			\rule{0pt}{12pt} & Error & Order & Error & Order & Error & Order\\
			\hline
			\rule{0pt}{12pt}$ 121$ & $6.4167\times 10^{-4}$ & - & $4.0826\times 10^{-4}$ & - & $5.0923\times 10^{-4}$ & -\\
			\rule{0pt}{12pt}$ 441$ & $1.6445\times 10^{-4}$ & $1.9823$ & $1.0243\times 10^{-4}$ & $2.0134$ & $1.3785\times 10^{-4}$ & $1.9190$ \\
			\rule{0pt}{12pt}$ 1681$ & $4.1914\times 10^{-5}$ & $1.9783$ & $2.5968\times 10^{-5}$ & $1.9876$ & N.C.* & -\\
			\hline
			\end{tabular}
	\label{P1 Error Chebyshev}
	\end{center}
\vspace{-0.3cm}
\scriptsize{*Not Convergent}
\end{table}

\begin{table}[h!]
		\caption{Comparison of CPU Time for Different Nodes per Star for Uniform Distribution of Nodes for Example 1}
		\vspace{-0.3 cm}
		\begin{center}
				\vspace{-0.4 cm}			
				\begin{tabular}{c c c c }
						\hline
						\rule{0pt}{12pt} Nodes & $N_G=5,~N_R=5$ &$N_G=9,~N_R=5$& $N_G=9,~N_R=9$\\
						\hline
						\rule{0pt}{12pt}$ 121$ & $2.5625\times 10^{-2}$ & $6.2500\times 10^{-2}$  & $7.8126\times 10^{-2}$\\
						\hline
						\rule{0pt}{12pt}$ 441$ & $0.2343$ &  $0.4062$ & $0.6093$\\
						\hline
						\rule{0pt}{12pt}$1861$ &$ 3.4687$ & $4.1718$ & $5.2812$ \\
						\hline
					\end{tabular}
				\label{P1 Time Uniform}
			\end{center}
	\end{table}

\begin{table}[h!]
		\caption{Comparison of CPU Time for Different Nodes per Star for Chebyshev Distribution of Nodes for Example 1}
		\vspace{-0.3 cm}
		\begin{center}
				\vspace{-0.4 cm}
				\begin{tabular}{c c c c}
						\hline
						\rule{0pt}{12pt} Nodes & $N_G=5,~N_R=5$ &$N_G=9,~N_R=5$& $N_G=9,~N_R=9$\\
						\hline
						\rule{0pt}{12pt}$ 121$ & $6.2500\times 10^{-2}$ & $9.3750\times 10^{-2}$ & $0.1562$\\
						\hline
						\rule{0pt}{12pt}$ 441$ & $2.6250$ &  $4.2968$ & $594.9375$ \\
						\hline
						\rule{0pt}{12pt}$1861$ &$ 115.2968$ & $144.75$ & - \\
						\hline
					\end{tabular}
				\label{P1 Time Chebyshev}
			\end{center}
	\end{table}

In order to study the influence of the number of nodes per star around the central node in the hybrid method, three different types of stars, namely, five nodes per star for both GFD and RBF ($N_G=5,~N_R=5$), nine nodes per star for GFD and five nodes for RBF ($N_G=9,~N_R=5$) and nine nodes per star for both GFD and RBF ($N_G=9,~N_R=9$) are considered. The RMS error and order of convergence of these three stars have been calculated for both uniform and Chebyshev distribution of nodes and are tabulated in Tables \ref{P1 Error Uniform} and \ref{P1 Error Chebyshev} respectively. For Chebyshev distribution of nodes, nodal spacing is calculated by using the formula:

\begin{equation}\nonumber
	h=\max_{1\leq i\leq n} \min_{1\leq j\leq n,~j\neq i}|\textbf x_i-\textbf x_j|.
\end{equation}
In case of uniform distribution, the error and order of convergence with ($N_G = 5,~N_R = 5$) and ($N_G = 9,~N_R = 5$) are comparable, showing no significant difference. On the other hand, ($N_G = 9,~N_R = 9$) lags behind the other two types of stars both in error and order of convergence. Similar results are noticed with Chebyshev distribution of nodes as can be seen in Table \ref{P1 Error Chebyshev}. The CPU time taken for convergence for all three types of stars are tabulated in Tables \ref{P1 Time Uniform} and \ref{P1 Time Chebyshev}. It is clear from both the tables that the CPU time taken by ($N_G = 5,~N_R = 5$) is much less than the other two stars, proving that ($N_G = 5,~N_R = 5$) produces the most efficient results.\\

\begin{table}[h!]
	\caption{Comparison of Error and Order of Convergence for Hybrid GFD-RBF and GFD for Uniform Distribution of Nodes for Example 1}
	\vspace{-0.3 cm}
	\begin{center}
		\vspace{-0.4 cm}
		\begin{tabular}{ccccc}
			\hline
			\rule{0pt}{12pt}	\multirow{2}{*}{Nodes} &\multicolumn{2}{c}{Hybrid GFD-RBF}  & \multicolumn{2}{c}{GFD} \\
			\cline{2-5}
			\rule{0pt}{12pt} & Error & Order & Error & Order\\
			\hline
			\rule{0pt}{12pt}$ 121$ & $3.4543\times 10^{-5}$ & - & $2.2081\times 10^{-4}$ & - \\
			\rule{0pt}{12pt}$ 441$ & $8.9838\times 10^{-6}$ & $1.9429$ & $5.7462\times 10^{-5}$  & $1.9421$\\
			\rule{0pt}{12pt}$ 1681$ & $2.2998\times 10^{-6}$ & $1.9658$ & $1.4693\times 10^{-5}$ & $1.9674$\\
			\hline
		\end{tabular}
		\label{P1 Error GFD Uniform}
	\end{center}
\end{table}

\begin{table}[h!]
	\caption{Comparison of Error and Order of Convergence for Hybrid GFD-RBF and GFD for Chebyshev Distribution of Nodes for Example 1}
	\vspace{-0.3 cm}
	\begin{center}
		\vspace{-0.4 cm}
		\begin{tabular}{ccccc}
			\hline
			\rule{0pt}{12pt}	\multirow{2}{*}{Nodes} &\multicolumn{2}{c}{Hybrid GFD-RBF}  & \multicolumn{2}{c}{GFD} \\
			\cline{2-5}
			\rule{0pt}{12pt} & Error & Order & Error & Order\\
			\hline
			\rule{0pt}{12pt}$ 121$ & $6.4167\times 10^{-4}$ & - & $6.5959\times 10^{-4}$ & - \\
			\rule{0pt}{12pt}$ 441$ & $1.6445\times 10^{-4}$ & $1.9823$ & $1.6942\times 10^{-4}$ & $1.9792$ \\
			\rule{0pt}{12pt}$ 1681$ & $4.19146\times 10^{-5}$ & $1.9783$ & $4.3183\times 10^{-5}$  & $1.9793$\\
			\hline
		\end{tabular}
		\label{P1 Error GFD Chebyshev}
	\end{center}
\end{table}

\begin{table}[h!]
	\caption{Comparison of CPU Time for Hybrid GFD-RBF and GFD for Uniform Distribution of Nodes for Example 1}
	\vspace{-0.3 cm}
	\begin{center}
		\vspace{-0.4 cm}
		\begin{tabular}{cccc}
			\hline
			\rule{0pt}{12pt}	\multirow{2}{*}{Nodes} & \multirow{2}{*}{Hybrid GFD-RBF}  & \multirow{2}{*}{GFD} & Time saved with Hybrid\\
			\rule{0pt}{12pt} &  & & (w.r.t. GFD)\\
			\hline
			\rule{0pt}{12pt}$ 121$ & $2.5625\times 10^{-2}$ & $3.1250\times 10^{-2}$ & $18\%$ \\
			\rule{0pt}{12pt}$ 441$ & $0.2343$ &  $0.2900$ & $19.20\%$ \\
			\rule{0pt}{12pt}$1861$  & $2.7656$  & $ 3.4687$ &$20.26\%$\\
			\hline
		\end{tabular}
		\label{P1 Time GFD Uniform}
	\end{center}
\end{table}

\begin{table}[h!]
	\caption{Comparison of CPU Time for Hybrid GFD-RBF and GFD for Chebyshev Distribution of Nodes for Example 1}
	\vspace{-0.3 cm}
	\begin{center}
		\vspace{-0.4 cm}
		\begin{tabular}{cccc}
			\hline
			\rule{0pt}{12pt}	\multirow{2}{*}{Nodes} & \multirow{2}{*}{Hybrid GFD-RBF}  & \multirow{2}{*}{GFD} & Time saved with Hybrid\\
			\rule{0pt}{12pt} &  & & (w.r.t. GFD)\\
			\hline
			\rule{0pt}{12pt}$ 121$ & $6.2500\times 10^{-2}$ & $9.3750\times 10^{-2}$ & $33.33\%$\\
			\rule{0pt}{12pt}$ 441$ & $2.6250$ &  $3.8593$ & $31.9824\%$  \\
			\rule{0pt}{12pt}$1861$ &$115.2968$ & $158.3750$  & $27.2001\%$\\
			\hline
		\end{tabular}
		\label{P1 Time GFD Chebyshev}
	\end{center}
\end{table}
The efficiency of the hybrid GFD-RBF ($N_G = 5,~N_R = 5$) is also compared with GFD ($N_G=9$) for both types of nodal distributions in the domain and tabulated in Tables \ref{P1 Error GFD Uniform}-\ref{P1 Time GFD Chebyshev}. It has been found from Tables \ref{P1 Error GFD Uniform} and \ref{P1 Error GFD Chebyshev} that hybrid method is producing more accurate results than GFD method. Tables \ref{P1 Time GFD Uniform} and \ref{P1 Time GFD Chebyshev} show that hybrid method is able to save a significant amount of CPU time compared with GFD method. The maximum CPU time saving with hybrid method is 20.26\% in the case of uniform distribution of nodes and 33.33\% in the case of Chebyshev nodes. Hence hybrid GFD-RBF $(N_G=5,N_R=5)$ method proves to be an efficient meshless method to solve the linear convection-diffusion problems.\\

\textbf{Example 2: Nonlinear Convection-Diffusion Equation}
$$\left(\frac{\partial U}{\partial x}\right)^2+\left(\frac{\partial U}{\partial y}\right)^2+U\left(\frac{\partial^2 U}{\partial x^2} +\frac{\partial^2 U}{\partial y^2}\right)=2U^4,~~~~~~~~0\leq x,y\leq 1.$$
Exact solution:
$$U(x,y)=\frac{1}{\sqrt{(x+1)^2+(y+1)^2}}.$$

\begin{table}[h!]
	\caption{Comparison of Error, Order of Convergence and Condition Number of the Stiffness Matrix for Example 2 with $\epsilon=0.6$}
	\vspace{-0.3 cm}
	\begin{center}
		\vspace{-0.4 cm}
		\begin{tabular}{c c c c c c c }
			\hline
			\rule{0pt}{12pt} \multirow{2}{*}{Nodes} & \multicolumn{3}{c}{Hybrid GFD-RBF}  &\multicolumn{3}{c}{CD2}\\
			\cline{2-7}
			\rule{0pt}{12pt} & Error & Order & C.N. & Error & Order & C.N. \\
			\hline
			\rule{0pt}{12pt} $121$ & $9.9818\times 10^{-6}$ & - & $45.1127$ & $5.1342\times 10^{-5}$ & - & $45.1101$ \\
			\rule{0pt}{12pt}$ 441$ & $2.6321\times 10^{-6}$ & $1.9320$ & $199.1313$ & $1.3633\times 10^{-5}$ & $1.9129$ & $199.1288$\\
			\rule{0pt}{12pt}$ 1681$ & $7.3073\times 10^{-7}$ & $1.8488$ & $852.0021$ & $3.59178\times10^{-6}$ & $1.9244$ & $851.9988$\\
			\hline
		\end{tabular}
		\label{P2 CD2}
	\end{center}
\end{table}

In Example 2, the hybrid GFD-RBF method has been applied to nonlinear convection-diffusion equation. In this example, the numerical solution obtained by the hybrid method ($N_G=5,~N_R=5$) is closer to the analytical solution than CD2 as observed in Table \ref{P2 CD2}. The order of convergence and condition number of the stiffness matrix are in agreement with each other for both the methods. The error produced by hybrid method is also in line with GFD method as mentioned in Reference\cite{GFD theory 5}.\\

\begin{table}[h!]
	\caption{Comparison of Error and Order of Convergence for Different Nodes per Star for Uniform Distribution of Nodes for Example 2 with $\epsilon=0.6$}
	\vspace{-0.3 cm}
	\begin{center}
		\vspace{-0.4 cm}
		\begin{tabular}{c c c c c c c}
			\hline
			\rule{0pt}{12pt}\multirow{2}{*}{Nodes} &\multicolumn{2}{c}{$N_G=5,~N_R=5$} &\multicolumn{2}{c}{$N_G=9,~N_R=5$}& \multicolumn{2}{c}{$N_G=9,~N_R=9$}\\
			\cline{2-7}
			\rule{0pt}{12pt} & Error & Order & Error & Order & Error & Order\\
			\hline
			\rule{0pt}{12pt}$ 121$ & $9.9818\times 10^{-6}$ & - & $4.5601\times 10^{-5}$ & - & $1.1396\times 10^{-4}$ & -\\
			\rule{0pt}{12pt}$ 441$ & $2.6321\times 10^{-6}$ & $1.9230$ & $1.1290\times 10^{-5}$ & $2.0142$ & $3.5334\times 10^{-5}$ & $1.6894$\\
			\rule{0pt}{12pt}$ 1681$ & $7.3073\times 10^{-7}$ & $1.8488$ & $1.9550\times 10^{-6}$ & $2.5298$ & $8.2526\times 10^{-6}$ & $2.0981$ \\
			\hline
		\end{tabular}
		\label{P2 Error Uniform}
	\end{center}
\end{table}

\begin{table}[h!]
	\caption{Comparison of Error and Order of Convergence for Different Nodes per Star for Chebyshev Distribution of Nodes for Example 2 with $\epsilon=0.5$}
	\vspace{-0.3 cm}
	\begin{center}
		\vspace{-0.4 cm}
		\begin{tabular}{c c c c c c c}
			\hline
			\rule{0pt}{12pt}\multirow{2}{*}{Nodes} &\multicolumn{2}{c}{$N_G=5,~N_R=5$} &\multicolumn{2}{c}{$N_G=9,~N_R=5$}& \multicolumn{2}{c}{$N_G=9,~N_R=9$}\\
			\cline{2-7}
			\rule{0pt}{12pt} & Error & Order & Error & Order & Error & Order\\
			\hline
			\rule{0pt}{12pt}$ 121$ & $2.7062\times 10^{-5}$ & - & $5.6769\times 10^{-5}$ & - & $5.8346\times 10^{-5}$ & -\\
			\rule{0pt}{12pt}$ 441$ & $6.5296\times 10^{-6}$ & $2.0702$ & $1.4599\times 10^{-5}$ & $1.9775$ & N.C. &- \\
			\rule{0pt}{12pt}$ 1681$ & $1.4256\times 10^{-6}$ & $2.2035$ & $3.8659\times 10^{-6}$ & $1.9240$ & N.C.  & - \\
			\hline
		\end{tabular}
		\label{P2 Error Chebyshev}
	\end{center}
\end{table}

\begin{table}[h!]
		\caption{Comparison of CPU Time for Different Nodes per Star for Uniform Distribution of Nodes for Example 2}
		\vspace{-0.3 cm}
		\begin{center}
				\vspace{-0.4 cm}
				\begin{tabular}{cccc}
						\hline
						\rule{0pt}{12pt} Nodes & $N_G=5,~N_R=5$ &$N_G=9,~N_R=5$& $N_G=9,~N_R=9$\\
						\hline
						\rule{0pt}{12pt}$ 121$ & $0.3281$ & $0.6875$ & $1.0937$ \\
						\hline
						\rule{0pt}{12pt}$ 441$ & $2.9687$ &  $6.8125$ & $11.2812$ \\
						\hline
						\rule{0pt}{12pt}$1861$ &$30.4375$ & $89.6562$ & $160.5156$  \\
						\hline
					\end{tabular}
				\label{P2 Time Uniform}
			\end{center}
	\end{table}

\begin{table}[h!]
		\caption{Comparison of CPU Time for Different Nodes per Star for Chebyshev Distribution of Nodes for Example 2}
		\vspace{-0.3 cm}
		\begin{center}
				\vspace{-0.4 cm}
				\begin{tabular}{cccc}
					\hline
						\rule{0pt}{12pt} Nodes & $N_G=5,~N_R=5$ &$N_G=9,~N_R=5$& $N_G=9,~N_R=9$\\
					\hline
						\rule{0pt}{12pt}$ 121$ & $0.4531$ & $1.7343$ & $1.9062$ \\
						\hline
						\rule{0pt}{12pt}$ 441$ & $14.0345$ &  $73.2812$ & -\\
						\hline
						\rule{0pt}{12pt}$1861$ &$ 684.2656$ & $4752.2656$ & - \\
						\hline
					\end{tabular}
				\label{P2 Time Chebyshev}
			\end{center}
	\end{table}

When the number of nodes in a star is increased, it is observed from Table \ref{P2 Error Uniform} that in case of uniform distribution, ($N_G=5,~N_R=5$) has less error than ($N_G=9,~N_R=5$) though the order of convergence is high in ($N_G=9,~N_R=5$). On the other hand, RMS error and order of convergence produced by ($N_G=9,~N_R=9$) are not as satisfactory as with the other two types of stars. In case of Chebyshev distribution of nodes, ($N_G=5,~N_R=5$) is giving more accurate results, as is seen from Table \ref{P2 Error Chebyshev}. It is also clear from Tables \ref{P2 Time Uniform} and \ref{P2 Time Chebyshev} that ($N_G=5,~N_R=5$) is computationally efficient as the CPU time taken is less than the other two.\\

\begin{table}[h!]
	\caption{Comparison of Error and Order of Convergence for Hybrid GFD-RBF and GFD for Uniform Distribution of Nodes for Example 2}
	\vspace{-0.3 cm}
	\begin{center}
		\vspace{-0.4 cm}
		\begin{tabular}{ccccc}
			\hline
			\rule{0pt}{12pt}	\multirow{2}{*}{Nodes} &\multicolumn{2}{c}{Hybrid GFD-RBF}  & \multicolumn{2}{c}{GFD} \\
			\cline{2-5}
			\rule{0pt}{12pt} & Error & Order & Error & Order\\
			\hline
			\rule{0pt}{12pt}$ 121$ & $9.9818\times 10^{-6}$ & - & $7.9030\times 10^{-5}$ & - \\
			\rule{0pt}{12pt}$ 441$ & $2.6321\times 10^{-6}$ & $1.9230$ & $2.0600\times 10^{-5}$ & $1.9397$ \\
			\rule{0pt}{12pt}$ 1681$ & $7.3073\times 10^{-7}$ & $1.8488$ & $5.1814\times 10^{-6}$ & $1.9952$ \\
			\hline
		\end{tabular}
		\label{P2 Error GFD Uniform}
	\end{center}
\end{table}

\begin{table}[h!]
	\caption{Comparison of Error and Order of Convergence for Hybrid GFD-RBF and GFD for Chebyshev Distribution of Nodes for Example 2}
	\vspace{-0.3 cm}
	\begin{center}
		\vspace{-0.4 cm}
		\begin{tabular}{ccccc}
			\hline
			\rule{0pt}{12pt}	\multirow{2}{*}{Nodes} &\multicolumn{2}{c}{Hybrid GFD-RBF}  & \multicolumn{2}{c}{GFD} \\
			\cline{2-5}
			\rule{0pt}{12pt} & Error & Order & Error & Order\\
			\hline
			\rule{0pt}{12pt}$ 121$ & $2.7062\times 10^{-5}$ & - & $9.8615\times 10^{-5}$ & - \\
			\rule{0pt}{12pt}$ 441$ & $6.5296\times 10^{-6}$ & $2.0702$ & $2.5439\times 10^{-5}$  & $1.9730$\\
			\rule{0pt}{12pt}$ 1681$ & $1.4256\times 10^{-6}$ & $2.2035$ & $5.8435\times 10^{-6}$  & $2.1299$\\
			\hline
		\end{tabular}
		\label{P2 Error GFD Chebyshev}
	\end{center}
\end{table}

\begin{table}[h!]
	\caption{Comparison of CPU Time for Hybrid GFD-RBF and GFD for Uniform Distribution of Nodes for Example 2}
	\vspace{-0.3 cm}
	\begin{center}
		\vspace{-0.4 cm}
		\begin{tabular}{cccc}
			\hline
			\rule{0pt}{12pt}	\multirow{2}{*}{Nodes} & \multirow{2}{*}{Hybrid GFD-RBF}  & \multirow{2}{*}{GFD}  & Time saved with Hybrid\\
			\rule{0pt}{12pt} &  & & (w.r.t. GFD)\\
			\hline
			\rule{0pt}{12pt}$ 121$ & $0.3281$ & $0.4218$ & $22.21\%$ \\
			\rule{0pt}{12pt}$ 441$ & $2.9687$ &  $3.8625$ & $23.14\%$ \\
			\rule{0pt}{12pt}$1861$ &$30.4375$ & $40.0312$ & $23.96\%$ \\
			\hline
		\end{tabular}
		\label{P2 Time GFD Uniform}
	\end{center}
\end{table}

\begin{table}[h!]
	\caption{Comparison of CPU Time for Hybrid GFD-RBF and GFD for Chebyshev Distribution of Nodes for Example 2}
	\vspace{-0.3 cm}
	\begin{center}
		\vspace{-0.4 cm}
		\begin{tabular}{cccc}
			\hline
			\rule{0pt}{12pt}	\multirow{2}{*}{Nodes} & \multirow{2}{*}{Hybrid GFD-RBF}  & \multirow{2}{*}{GFD}  & Time saved with Hybrid\\
			\rule{0pt}{12pt} &  & & (w.r.t. GFD)\\
			\hline
			\rule{0pt}{12pt}$ 121$ & $0.4531$ & $0.8437$ & $46.29\%$ \\
			\rule{0pt}{12pt}$ 441$ & $14.0345$ &  $26.7343$ & $47.50\%$ \\
			\rule{0pt}{12pt}$1861$ &$684.2656$ & $1329.9220$ & $48.54\%$ \\
			\hline
		\end{tabular}
		\label{P2 Time GFD Chebyshev}
	\end{center}
\end{table}

The RMS error and the order of convergence of hybrid method ($N_G=5,~N_R=5$) are compared with GFD for both types of nodal distributions and results are tabulated in Tables \ref{P2 Error GFD Uniform} and \ref{P2 Error GFD Chebyshev}. It is found from the tables that the hybrid method gives more accurate results than GFD method while maintaining the same level of order of convergence. Comparison of CPU time for both the methods is also done which is given in Tables \ref{P2 Time GFD Uniform} and \ref{P2 Time GFD Chebyshev}. It is found that hybrid method is able to save a maximum CPU time of 23.96\% in the case of uniform nodes and 48.54\% in the case of Chebyshev nodes. Hence, the developed hybrid GFD-RBF $(N_G=5,N_R=5)$ method solves nonlinear convection-diffusion problem efficiently.\\

\textbf{Example 3: Coupled Nonlinear Convection-Diffusion Equation}
\begin{equation}\nonumber
	\left.
	\begin{aligned}
		U\frac{\partial U}{\partial x}+V\frac{\partial U}{\partial y}-\frac{\partial^2 U}{\partial x^2} -\frac{\partial^2 U}{\partial y^2}&=f_1(x,y), \\
		U\frac{\partial V}{\partial x}+V\frac{\partial V}{\partial y}-\frac{\partial^2 V}{\partial x^2} -\frac{\partial^2 V}{\partial y^2}&=f_2(x,y),
	\end{aligned}
	\right.
	\:\:\:\:\: \:\:\:0\leq x,y\leq \pi.
\end{equation}
Exact solution:
$$U(x,y)=-\cos x\sin y,~~~~V(x,y)=\sin x\cos y.$$

\begin{table}[h!]
	\caption{Comparison of Error, Order of Convergence and Condition Number of the Stiffness Matrix for Example 3 with $\epsilon=0.3$}
	\vspace{-0.3 cm}
	\begin{center}
		\vspace{-0.4 cm}
		\scriptsize
		\begin{tabular}{p{0.05cm} c c c c c c c c c}
		
			\hline
			\rule{0pt}{12pt} & \multirow{2}{*}{Nodes} & \multicolumn{3}{c}{Hybrid GFD-RBF}  &\multicolumn{3}{c}{CD2} &\multicolumn{2}{c}{Ref.\cite{Comparision}(RBF-FD)}\\
			\cline{3-10}
			\rule{0pt}{12pt} & & Error & Order & C.N. & Error & Order & C.N. & Error & Order \\
			\hline
			\rule{0pt}{12pt}\multirow{3}{*}{$u$} & $121$ & $3.6074\times 10^{-4}$ & - & $38.8377$ & $9.7989\times 10^{-4}$ & - & $38.8193$ & $4.8366\times 10^{-4}$ & -\\
			\rule{0pt}{12pt} & $ 441$ & $9.2307\times 10^{-5}$ & $1.9664$ & $156.8412$ & $2.5603\times 10^{-4}$ & $1.9363$ & $156.8192$ & $1.1976\times 10^{-4}$ & $2.01$\\
			\rule{0pt}{12pt} & $ 1681$ & $2.3507\times 10^{-5}$ & $1.9733$ & $628.8272$ & $6.5524\times10^{-5}$ & $1.9662$ & $628.8043$ & $3.0248\times 10^{-5}$ & $2.71$\\
			\hline
			\rule{0pt}{12pt}\multirow{3}{*}{$v$} & $121$ & $6.5152\times 10^{-4}$ & - & $38.8377$ & $2.5574\times 10^{-4}$ & - & $38.8193$ & $1.0350\times 10^{-3}$ & - \\
			\rule{0pt}{12pt} & $ 441$ & $1.7188\times 10^{-4}$ & $1.9224$ & $156.8412$ & $6.8869\times 10^{-5}$ & $1.8927$ & $156.8192$ & $ 2.7923\times 10^{-4}$ & $1.89$\\
			\rule{0pt}{12pt} & $ 1681$ & $4.4140\times 10^{-5}$ & $1.9612$ & $628.8272$ & $1.7715\times10^{-5}$ & $1.9588$ & $628.8043$ & $2.4539\times 10^{-4}$ & $0.19$\\
			\hline	
		\end{tabular}
		\label{P3 CD2}
	\end{center}
\end{table}

The developed hybrid method is also tested on coupled nonlinear convection-diffusion equation in Example 3. The comparison of the hybrid method is done with CD2 and RBF-FD\cite{Comparision} methods in terms of error, order of convergence and condition number of the stiffness matrix and the results are tabulated in Table \ref{P3 CD2}. It is found the numerical solution obtained by the proposed hybrid method is more accurate that the other two methods. The order of convergence and condition number of the hybrid method is in agreement with the other two methods.\\

\begin{table}[h!]
	\caption{Comparison of Error and Order of Convergence for Different Nodes per Star for Uniform Distribution of Nodes for Example 3 with $\epsilon=0.3$}
	\vspace{-0.3 cm}
	\begin{center}
		\vspace{-0.4 cm}
		\begin{tabular}{c c c c c c c c}
			\hline
			\rule{0pt}{12pt} &\multirow{2}{*}{Nodes} &\multicolumn{2}{c}{$N_G=5,~N_R=5$} &\multicolumn{2}{c}{$N_G=9,~N_R=5$}& \multicolumn{2}{c}{$N_G=9,~N_R=9$}\\
			\cline{3-8}
			\rule{0pt}{12pt}&  & Error & Order & Error & Order & Error & Order\\
			\hline
			\rule{0pt}{12pt} \multirow{3}{*}{$u$} & $121$ & $3.6074\times 10^{-4}$ & - &  $1.3650\times 10^{-3}$ & - &  $1.1579\times 10^{-2}$ & - \\
			\rule{0pt}{12pt} & $441$ & $9.2307\times 10^{-5}$ & $1.9664$ & $3.5707\times 10^{-4}$ & $1.9346$ & $5.4368\times 10^{-3}$ & $1.0906$  \\
			\rule{0pt}{12pt} & $1681$ & $2.3507\times 10^{-5}$ & $1.9733$ & $9.1420\times 10^{-5}$ & $1.9656$ & $2.6093\times 10^{-3}$ & $1.0590$\\
			\hline
			\rule{0pt}{12pt} \multirow{3}{*}{$v$} & $121$ & $6.5152\times 10^{-4}$ & - &  $3.7985\times 10^{-4}$& -  & $9.8660\times 10^{-3}$ & -\\
			\rule{0pt}{12pt} & $441$ &  $1.7188\times 10^{-4}$ & $1.9224$  & $9.7175\times 10^{-5}$ & $1.9667$ & $3.2541\times 10^{-3}$ & $1.6002$ \\
			\rule{0pt}{12pt} & $1681$ & $4.4140\times 10^{-5}$ & $1.9612$ & $2.4759\times 10^{-5}$ & $1.9726$ &  $1.0056\times 10^{-3}$ & $1.6942$\\
			\hline
		\end{tabular}
		\label{P3 Error Uniform}
	\end{center}
\end{table}

\begin{table}[h!]
	\caption{Comparison of Error and Order of Convergence for Different Nodes per Star for Chebyshev Distribution of Nodes for Example 3 with $\epsilon=0.3$}
	\vspace{-0.3 cm}
	\begin{center}
		\vspace{-0.4 cm}
		\begin{tabular}{c c c c c c c c}
			\hline
			\rule{0pt}{12pt} &\multirow{2}{*}{Nodes} &\multicolumn{2}{c}{$N_G=5,~N_R=5$} &\multicolumn{2}{c}{$N_G=9,~N_R=5$}& \multicolumn{2}{c}{$N_G=9,~N_R=9$}\\
			\cline{3-8}
			\rule{0pt}{12pt}&  & Error & Order & Error & Order & Error & Order\\
			\hline
			\rule{0pt}{12pt} \multirow{3}{*}{$u$} & $121$ & $8.4971\times 10^{-4}$ & - & $1.8755\times 10^{-3}$ & - & $1.5195\times 10^{-2}$ & -\\
			\rule{0pt}{12pt} & $441$ & $2.1565\times 10^{-4}$ & $1.9967$ & $4.8936\times 10^{-4}$ & $1.9564$ & N.C. & - \\
			\rule{0pt}{12pt} & $1681$ & $5.4810\times 10^{-5}$ & $1.9834$ & $1.2522\times 10^{-4}$ & $1.9737$ & N.C. & - \\
			\hline
			
			\rule{0pt}{12pt} \multirow{3}{*}{$v$} & $121$ & $9.8395\times 10^{-4}$ & - & $4.8010\times 10^{-4}$& - & $2.6727\times 10^{-2}$ & -\\
			\rule{0pt}{12pt} & $441$ & $2.8839\times 10^{-4}$ & $1.7871$ & $1.2765\times 10^{-4}$& $1.8954$ & N.C. & - \\
			\rule{0pt}{12pt} & $1681$ & $7.6075\times 10^{-5}$ & $1.9283$ & $3.2857\times 10^{-5}$ & $1.9651$ & N.C. & - \\
			\hline
		\end{tabular}
		\label{P3 Error Chebyshev}
	\end{center}
\end{table}

\begin{table}[h!]
		\caption{Comparison of CPU Time for Different Nodes per Star for Uniform Distribution of Nodes for Example 3}
		\vspace{-0.3 cm}
		\begin{center}
				\vspace{-0.4 cm}
				\begin{tabular}{cccc}
						\hline
						\rule{0pt}{12pt} Nodes & $N_G=5,~N_R=5$ & $N_G=9,~N_R=5$& $N_G=9,~N_R=9$\\
						\hline
						\rule{0pt}{12pt}$ 121$ & $3.1250\times 10^{-2}$ & $0.1562$ & $0.4062$ \\
						\hline
						\rule{0pt}{12pt}$ 441$ & $0.5781$ &  $2.2500$ & $7.0156$\\
						\hline
						\rule{0pt}{12pt}$1861$ &$7.8281$ & $35.2500$ & $104.2500$ \\
						\hline
					\end{tabular}
				\label{P3 Time Uniform}
			\end{center}
	\end{table}

\begin{table}[h!]
		\caption{Comparison of CPU Time for Different Nodes per Star for Chebyshev Distribution of Nodes for Example 3}
		\vspace{-0.3 cm}
		\begin{center}
				\vspace{-0.4 cm}
				\begin{tabular}{cccc}
						\hline
						\rule{0pt}{12pt} Nodes & $N_G=5,~N_R=5$ & $N_G=9,~N_R=5$& $N_G=9,~N_R=9$\\
						\hline
						\rule{0pt}{12pt}$ 121$ & $0.1718$ & $0.4843$ & $0.8593$ \\
						\hline
						\rule{0pt}{12pt}$ 441$ & $5.8593$ &  $22.0156$ & -\\
						\hline
						\rule{0pt}{12pt}$1861$ &$ 271.2187$ & $1396.4687$ & - \\
						\hline
					\end{tabular}
				\label{P3 Time Chebyshev}
			\end{center}
	\end{table}

The influence of the number of nodes per star in the hybrid method has been studied for uniform and Chebyshev distribution of nodes in the domain. Table \ref{P3 Error Uniform} shows that, ($N_G=5,~N_R=5$) has less error than ($N_G=9,~N_R=5$) in case of uniform distribution  of nodes while maintaining the same level of order of convergence. On the other hand, ($N_G=9,~N_R=9$) has low performance than the other two methods. Similar results are there for Chebyshev distribution of nodes, as is seen from Table \ref{P3 Error Chebyshev}. The CPU time taken by ($N_G=5,~N_R=5$) is also less than the other two stars in both types of nodal distributions, as observed from Tables \ref{P3 Time Uniform} and \ref{P3 Time Chebyshev}, proving its efficiency over the other two methods.\\

\begin{table}[h!]
	\caption{Comparison of Error and Order of Convergence for Hybrid GFD-RBF and GFD for Uniform Distribution of Nodes for Example 3}
	\vspace{-0.3 cm}
	\begin{center}
		\vspace{-0.4 cm}
		\begin{tabular}{cccccc}
			\hline
			\rule{0pt}{12pt}	&\multirow{2}{*}{Nodes} &\multicolumn{2}{c}{Hybrid GFD-RBF}  &\multicolumn{2}{c}{GFD}\\
			\cline{3-6}
			\rule{0pt}{12pt} & & Error & Order & Error & order\\
			\hline
			\rule{0pt}{12pt}\multirow{3}{*}{$u$} & $121$ & $3.6074\times 10^{-4}$ & - & $2.2484\times 10^{-3}$ & - \\
			\rule{0pt}{12pt} & $441$ & $9.2307\times 10^{-5}$ & $1.9664$ & $5.8221\times 10^{-4}$ & $1.9492$ \\
			\rule{0pt}{12pt} & $1681$ & $2.3507\times 10^{-5}$ & $1.9733$ & $1.4865\times 10^{-4}$ & $1.9676$ \\
			\hline
			\rule{0pt}{12pt}\multirow{3}{*}{$v$} & $121$  & $6.5152\times 10^{-4}$ & - & $8.7300\times 10^{-4}$ & -\\
			\rule{0pt}{12pt} & $441$ & $1.7188\times 10^{-4}$ & $1.9224$ & $2.3464\times 10^{-4}$ & $1.8955$ \\
			\rule{0pt}{12pt} & $1681$  & $4.4140\times 10^{-5}$ & $1.9612$ & $6.0088\times 10^{-5}$ & $1.9653$\\
			\hline
		\end{tabular}
		\label{P3 Error GFD Uniform}
	\end{center}
\end{table}

\begin{table}[h!]
	\caption{Comparison of Error and Order of Convergence for Hybrid GFD-RBF and GFD for Chebyshev Distribution of Nodes for Example 3}
	\vspace{-0.3 cm}
	\begin{center}
		\vspace{-0.4 cm}
		\begin{tabular}{cccccc}
			\hline
			\rule{0pt}{12pt}	&\multirow{2}{*}{Nodes} &\multicolumn{2}{c}{Hybrid GFD-RBF}  &\multicolumn{2}{c}{GFD}\\
			\cline{3-6}
			\rule{0pt}{12pt} & & Error & Order & Error & order\\
			\hline
			\rule{0pt}{12pt} \multirow{3}{*}{$u$} & $121$ & $8.4971\times 10^{-4}$ & - & $3.1552\times 10^{-3}$ & - \\
			\rule{0pt}{12pt} & $441$ & $2.1565\times 10^{-4}$ & $1.9967$ & $8.0597\times 10^{-4}$ & $1.9873$ \\
			\rule{0pt}{12pt} & $1681$ & $5.4810\times 10^{-5}$ & $1.9834$ & $2.0513\times 10^{-4}$ & $1.9769$ \\
			\hline
			
			\rule{0pt}{12pt} \multirow{3}{*}{$v$} & $121$ & $8.4971\times 10^{-4}$ & - & $1.6287\times 10^{-3}$ & - \\
			\rule{0pt}{12pt} & $441$ &  $2.8839\times 10^{-4}$ & $1.7871$ & $4.1602\times 10^{-4}$ & $1.9873$ \\
			\rule{0pt}{12pt} & $1681$  & $7.6075\times 10^{-5}$ & $1.9283$ & $1.0584\times 10^{-4}$ & $1.9822$ \\
			\hline
		\end{tabular}
		\label{P3 Error GFD Chebyshev}
	\end{center}
\end{table}

\begin{table}[h!]
	\caption{Comparison of CPU Time for Hybrid GFD-RBF and GFD for Uniform Distribution of Nodes for Example 3}
	\vspace{-0.3 cm}
	\begin{center}
		\vspace{-0.4 cm}
		\begin{tabular}{cccc}
			\hline
			\rule{0pt}{12pt}	\multirow{2}{*}{Nodes} & \multirow{2}{*}{Hybrid GFD-RBF}  & \multirow{2}{*}{GFD}  & Time saved with Hybrid\\
			\rule{0pt}{12pt} &  & & (w.r.t. GFD)\\
			\hline
			\rule{0pt}{12pt}$ 121$ & $3.1250\times10^{-2}$ & $3.6825\times10^{-2}$ & $15.13\%$\\
			\rule{0pt}{12pt}$ 441$ & $0.5781$ &  $0.5312$ & $-8.82\%$  \\
			\rule{0pt}{12pt}$1861$ &$7.8281$ & $7.1187$ & $-9.96\%$ \\
			\hline
		\end{tabular}
		\label{P3 Time GFD Uniform}
	\end{center}
\end{table}

\begin{table}[h!]
	\caption{Comparison of CPU Time for Hybrid GFD-RBF and GFD for Chebyshev Distribution of Nodes for Example 3}
	\vspace{-0.3 cm}
	\begin{center}
		\vspace{-0.4 cm}
		\begin{tabular}{cccc}
			\hline
			\rule{0pt}{12pt}	\multirow{2}{*}{Nodes} & \multirow{2}{*}{Hybrid GFD-RBF}  & \multirow{2}{*}{GFD}  & Time saved with Hybrid\\
			\rule{0pt}{12pt} &  & & (w.r.t. GFD)\\
			\hline
			\rule{0pt}{12pt}$ 121$ & $0.1718$ & $0.2187$ & $21.44\%$\\
			\rule{0pt}{12pt}$ 441$ & $5.8593$ &  $9.7187$ & $39.71\%$ \\
			\rule{0pt}{12pt}$1861$ &$271.2187$ & $516.0625$ & $47.44\%$ \\
			\hline
		\end{tabular}
		\label{P3 Time GFD Chebyshev}
	\end{center}
\end{table}
Numerical solution obtained by hybrid method ($N_G=5,~N_R=5$) is more accurate than GFD as is evident from Tables \ref{P3 Error GFD Uniform} and \ref{P3 Error GFD Chebyshev}. With the increase of number of nodes in the uniform distribution, it is observed that the hybrid method needs more CPU time than GFD as shown in Table \ref{P3 Time GFD Uniform}. But it is able to save a maximum of 47.44\% CPU time in case of Chebyshev distribution of nodes as mentioned in Table \ref{P3 Time GFD Chebyshev}. So in general, for solving coupled nonlinear convection-diffusion problems, hybrid GFD-RBF $(N_G=5,N_R=5)$ proved to be an efficient meshless method.\\

\section{Conclusion}
Hybrid GFD-RBF method has been tested on linear and nonlinear convection-diffusion problems for both uniform and non-uniform distributions of nodes. For uniform node distribution, it was compared with the second-order central finite difference (CD2) method. Results show that the hybrid method produces significantly smaller errors than CD2 across all tested problems. Furthermore, the condition number of the stiffness matrix obtained from the hybrid method was evaluated and found to be satisfactory in comparison to CD2. It was observed that the hybrid method achieves optimal results with five nodes per star $(N_G=5,~N_R=5)$ for both uniform and non-uniform node distributions, whereas the GFD method requires nine nodes. Despite using fewer nodes per star, the hybrid method gives more accurate results than the GFD method. Additionally, in cases involving non-uniform node distributions, the hybrid method requires less CPU time than the GFD method. The hybrid method was determined to be numerically second order accurate. But the limitation of the method is that it is not helpful if the PDE is of homogenous type. \\
So, the hybrid method would be more effective if the first order derivatives (convective terms) are discretized by GFD and higher order derivatives (diffusive and other higher order terms) by RBF-FD to solve non-homogenous PDEs. Hence the developed method is computationally efficient and hence useful in solving convection-diffusion equations especially where the convective terms are highly nonlinear.


\begin{thebibliography}{1}

	\bibitem{1}
	Hongwei, G., Shan. L., Hong, Z. [2023] \enquote{GMLS-based numerical manifold method in mechanical analysis of thin plates with complicated shape or cutouts,} \textit{Eng. Anal. Bound. Elem.} \textbf{151}, 597-623.
	
	\bibitem{2}
	Xitailang, C., Shan, L., Zenglong, L., Hongwei, G.,	Hong, Z. [2024] \enquote{Meshless numerical manifold method with novel subspace tracking and CSS locating techniques for slope stability analysis,} \textit{Comput. Geotech.} \textbf{166}, 106025.
	
	\bibitem{GFD theory 1}
	Jensen, P.S. [1972] \enquote{Finite difference technique for variable grid,} \textit{Comput. \& amp; Structures} \textbf{2}, 17–29.
 \bibitem{GFD theory 2}
	Perrone, N., Kao, R. [1975] \enquote{A general finite difference method for arbitrary meshes,} \textit{Comput. \& amp; Structures} \textbf{5}, 45–58.

	\bibitem{GFD theory 3}
	Benito, J.J.,  Ureña, F., Gavete, L. [2001] \enquote{Influence several factors in the generalized finite difference method,} \textit{Appl. Math. Modelling} \textbf{25}, 1039–1053.
%
	\bibitem{GFD theory 4}
	Benito, J.J., Ureña, F., Gavete, L., Alonso, B. [2009] \enquote{Application of the generalized finite difference method to improve the approximated solution of pdes,} \textit{Comput. Model. Eng. Sci.} \textbf{38}, 39–58.
%
    \bibitem{GFD theory 6}
	Zheng, Z., Li, X. [2022] \enquote{Theoretical analysis of the generalized finite difference method,} \textit{Comput. Math. Appl.} \textbf{120}, 1-14.

	\bibitem{GFD app 2}
	Mochnacki, B., Majchrzak, E. [2010] \enquote{Numerical modeling of casting solidification using generalized finite difference method,} \textit{Mater. Sci. Forum} \textbf{638–642}, 2676–2681.
%
	\bibitem{GFD app 3}
	Ureña, F., Benito, J.J., Gavete, L. [2011] \enquote{Application of the generalized finite difference method to solve the advection-diffusion equation,} \textit{J. Comput. Appl. Math.} \textbf{235}, 1849–1855.
	
		\bibitem{GFD theory 5}
	Gavete, L., Ureña, F., Benito, J.J., García, A., Ureña, M., Salete, E. [2017] \enquote{Solving second order non-linear elliptic partial differential equations using generalized finite difference method,} \textit{J. Comput. Appl. Math.} \textbf{318}, 378–387.
	
	\bibitem{para}
	Ureña, F., Gavete, L., García, A., Benito, J.J., Vargas, A.M. [2019] \enquote{Solving second order non-linear parabolic PDEs using generalized finite difference method (GFDM),} \textit{J. Comput. Appl. Math.} \textbf{354}, 221-241.
	
	 
	\bibitem{hyper}
	 Ureña, F., Gavete, L., García, A., Benito, J.J., Vargas, A.M. [2020] \enquote{Solving second order non-linear hyperbolic PDEs using generalized finite difference method (GFDM),} \textit{J. Comput. Appl. Math.} \textbf{363}, 1-21.
	
	\bibitem{GFD app 4}
	Haung, T., Zhao, H., Chen, H., Yao, Y., Yu, P. [2022] \enquote{A hybrid cartesian-meshless method for the simulatrion of thermal flows with complex immersed objects,} \textit{Phys. Fluids} \textbf{34}, 103318.

	\bibitem{GFD app 5}
	Wright, G.B., Jones, A., Shankar, V. [2023] \enquote{MGM: A meshfree geometric multilevel method for systems arising from elliptic equations on point cloud surfaces,} \textit{SIAM J. Sci. Comput.} \textbf{45}(2), A312-A337.

	\bibitem{RBF theory 1}
	Kansa, E.J. [1990] \enquote{Multiquadrics—A scattered data approximation scheme with applications to computational fluid dynamics—II solutions to parabolic, hyperbolic and elliptic partial differential equations,} \textit{Comput. Math. Appl.} \textbf{19}(8–9), 147-161.

   \bibitem{RBF theory 2}
	Wright, G.B., Fornberg, B. [2006] \enquote{Scattered node compact finite difference-type formulas generated from radial basis functions,} \textit{J. Comput. Phys.} \textbf{212}(1), 99–123.

	\bibitem{RBF theory 3}
	Chandhini, G., Sanyasiraju, Y.V.S.S. [2007] \enquote{Local RBF-FD solutions for steady convection–diffusion problems,} \textit{Internat. J. Numer. Methods Engrg.} \textbf{72}, 352–378.
	
	\bibitem{RBF theory 4}
	Fornberg, B., Flyer, N. [2015] \enquote{Solving PDEs with radial basis functions,} \textit{Acta Numer.} \textbf{24}, 215–258.
	
 \bibitem{RBF app 1}
	Gunderman, D., Flyer, N., Fornberg, B. [2022] \enquote{Transport schemes in spherical geometries using spline-based RBF-FD with polynomials,} \textit{J. Comput. Phys.} \textbf{408}, 109256.
	
	\bibitem{RBF app 2}
	Shahane, S., Radhakrishnan, A., Vanka, S.P. [2021] \enquote{A high-order accurate meshless method for solution of incompressible fluid flow problems,} \textit{J. Comput. Phys.} \textbf{445}, 110623.
	
	\bibitem{RBF app 3}
	Sanyasiraju, Y.V.S.S., Chandhini, G. [2009]  \enquote{A note on two upwind strategies for RBF-based grid-free schemes to solve steady convection–diffusion equations,} \textit{Internat. J. Numer. Methods Fluids} \textbf{61}, 1053–1062.

	\bibitem{RBF app 4}
	 N.B. Barik, T.V.S. Sekhar, An Efficient Local RBF Meshless Scheme for Steady Convection–Diﬀusion Problems, Int. J. Comput. Methods 14 (2017) 6 1750064.
	
	\bibitem{RBF app 5}
	Barik, N.B., Sekhar, T.V.S. [2021] \enquote{Mesh-free multilevel iterative algorithm for Navier–Stokes equations,} \textit{Numer. Heat Transf. B: Fundam.} \textbf{79}(3), 150-164. 

	\bibitem{selection of points}
	Liszka, T., Orkisz, J. [1980] \enquote{The finite difference method at arbitrary irregular grids and its application in applied mechanics,} \textit{Comput. \& amp; Structures} \textbf{11}, 83-95.

	\bibitem{Rbf extra}
	Micchelli, C.A. [1986] \enquote{Interpolation of scattered data: distance matrices and conditionally positive definite functions,} \textit{Const. Approx.} \textbf{2}, 11–22.

	\bibitem{Multiquadric}
	Franke, R. [1982] \enquote{Scattered data interpolation: tests of some methods,} \textit{Math. Comput.} \textbf{38}, 181–200.
	
	\bibitem{Comparision}
	Barik, N.B., Sekhar, T.V.S. [2021] \enquote{A modified multilevel meshfree algorithm for steady convection-diffusion problems,} \textit{Int. J. Numer. Meth. Fluids.} \textbf{93}(7), 2121-2135.
\end{thebibliography}

\end{document}